\documentclass[12pt,a4paper]{article}
\usepackage{amsmath}
\usepackage{amssymb}
\usepackage{t1enc}
\usepackage[latin1]{inputenc}
\usepackage[english]{babel}
\usepackage{amsfonts}
\usepackage{latexsym}
\usepackage{bm}
\newtheorem{theorem}{Theorem}[section]

\newtheorem{conj}[theorem]{Conjecture}

\begin{document}
\title{Stars on trees}

\author{Peter Borg\\[5mm]
Department of Mathematics \\
University of Malta\\
Malta\\
\texttt{peter.borg@um.edu.mt} 
}

\date{}
\maketitle

\begin{abstract}
For a positive integer $r$ and a vertex $v$ of a graph $G$, let $\mathcal{I}_G^{(r)}(v)$ denote the set of all independent sets of $G$ that have exactly $r$ elements and contain $v$. 
Hurlbert and Kamat conjectured that for any $r$ and any tree $T$, there exists a leaf $z$ of $T$ such that $|\mathcal{I}_T^{(r)}(v)| \leq |\mathcal{I}_T^{(r)}(z)|$ for each vertex $v$ of $T$. They proved the conjecture for $r \leq 4$. For any $k \geq 3$, we construct a tree $T_k$ that has a vertex $x$ such that $x$ is not a leaf of $T_k$, $|\mathcal{I}_{T_k}^{(r)}(z)| < |\mathcal{I}_{T_k}^{(r)}(x)|$ for any leaf $z$ of $T_k$ and any $5 \leq r \leq 2k+1$, and $2k+1$ is the largest integer $s$ for which $\mathcal{I}_{T_k}^{(s)}(x)$ is non-empty. 
Therefore, the conjecture is not true for $r \geq 5$. 
\end{abstract}

\section{Introduction}

We shall use small letters such as $x$ to denote non-negative integers or elements of a set, capital letters such as $X$ to denote sets or graphs, and calligraphic letters such as $\mathcal{F}$ to denote \emph{families} (that is, sets whose members are sets themselves). The set $\{1, 2, \dots\}$ of positive integers is denoted by $\mathbb{N}$. For any $m, n \in \mathbb{N}$, the set $\{i \in \mathbb{N} \colon m \leq i \leq n\}$ is denoted by $[m,n]$, and we abbreviate $[1,n]$ to $[n]$. For a set $X$, the family $\{A \subseteq X \colon |A| = r\}$ of all $r$-element subsets of $X$ is denoted by ${X \choose r}$. 
%
If $x \in X$ and $\mathcal{F}$ is a family of subsets of $X$, then the family $\{F \in \mathcal{F} \colon x \in F\}$ is denoted by $\mathcal{F}(x)$ and is called a \emph{star of $\mathcal{F}$}. 
All arbitrary sets are assumed to be finite. 

A \emph{graph} $G$ is a pair $(X,Y)$, where $X$ is a set, called the \emph{vertex set of $G$}, and $Y$ is a subset of ${X \choose 2}$ and is called the \emph{edge set of $G$}. The vertex set of $G$ and the edge set of $G$ are denoted by $V(G)$ and $E(G)$, respectively. An element of $V(G)$ is called a \emph{vertex of $G$}, and an element of $E(G)$ is called an \emph{edge of $G$}. We may represent an edge $\{v,w\}$ by $vw$. 
A vertex $v$ of $G$ is a \emph{leaf of $G$} if there exists exactly one vertex $w$ of $G$ such that $vw \in E(G)$. 


If $H$ is a graph such that $V(H) \subseteq V(G)$ and $E(H) \subseteq E(G)$, then 
we say that \emph{$G$ contains $H$}. 

If $n \geq 2$ and $v_1, v_2, \dots, v_n$ are the distinct vertices of a graph $G$ with $E(G) = \{v_iv_{i+1} \colon i \in [n-1]\}$, then $G$ is called a \emph{$(v_1,v_n)$-path} or simply a \emph{path}. 

%
%



A graph $G$ is a \emph{tree} if $|V(G)| \geq 2$ and $G$ contains exactly one $(v,w)$-path for every $v, w \in V(G)$ with $v \neq w$. 


Let $G$ be a graph. A subset $I$ of $V(G)$ is an \emph{independent set of $G$} if $vw \notin E(G)$ for every $v, w \in I$. Let $\mathcal{I}_G^{(r)}$ denote the family of all independent sets of $G$ of size $r$. An independent set $J$ of $G$ is \emph{maximal} if $J \nsubseteq I$ for each independent set $I$ of $G$ such that $I \neq J$. The size of a smallest maximal independent set of $G$ is denoted by $\mu(G)$.

Hurlbert and Kamat \cite{HK} conjectured that for any $r \geq 1$ and any tree $T$, there exists a leaf $z$ of $T$ such that $\mathcal{I}_T^{(r)}(z)$ is a star of $\mathcal{I}_T^{(r)}$ of maximum size.

\begin{conj}[{\cite[Conjecture~1.25]{HK}}] \label{conj} For any $r \geq 1$ and any tree $T$, there exists a leaf $z$ of $T$ such that $|\mathcal{I}_T^{(r)}(v)| \leq |\mathcal{I}_T^{(r)}(z)|$ for each $v \in V(T)$.
\end{conj}
Hurlbert and Kamat \cite{HK} also showed that the conjecture is true for $r \leq 4$. In the next section, we show that the conjecture is not true for $r \geq 5$. For any $k \geq 3$, we construct a tree $T_k$ that has a vertex $x$ such that $x$ is not a leaf of $T_k$, $|\mathcal{I}_{T_k}^{(r)}(z)| < |\mathcal{I}_{T_k}^{(r)}(x)|$ for any leaf $z$ of $T_k$ and any $r \in [5,2k+1]$, and $2k+1$ is the largest integer $s$ for which $\mathcal{I}_{T_k}^{(s)}(x)$ is non-empty.

Conjecture~\ref{conj} was motivated by a problem of Holroyd and Talbot \cite{Holroyd,HT}. A family $\mathcal{A}$ is \emph{intersecting} if every two sets in $\mathcal{A}$ intersect.  We say that $\mathcal{I}_G^{(r)}$ has the \emph{star property} if at least one of the largest intersecting subfamilies of $\mathcal{I}_G^{(r)}$ is a star of $\mathcal{I}_G^{(r)}$.  Holroyd and Talbot introduced the problem of determining whether $\mathcal{I}_G^{(r)}$ has the star property for a given graph $G$ and an integer $r \geq 1$. 
The Holroyd--Talbot (HT) Conjecture \cite[Conjecture~7]{HT} claims that $\mathcal{I}_G^{(r)}$ has the star property if $\mu(G) \geq 2r$. 
By the classical Erd\H{o}s--Ko--Rado Theorem \cite{EKR}, the HT Conjecture is true if $G$ has no edges. The HT Conjecture has been verified for certain graphs \cite{BH1,BH,HST,HT,HK,T,W}. It is also verified in \cite{Borg} for any graph $G$ with $\mu(G)$ sufficiently large depending on $r$; this is the only result known for the case where $G$ is a tree that is not a path (the problem for paths is solved in \cite{HST}), apart from the fact that $\mathcal{I}_G^{(r)}$ may not have the star property for certain values of $r$ (indeed, if $G$ is the tree $(\{0\} \cup [n], \{\{0,i\} \colon i \in [n]\})$ and $2 \leq n/2 < r < n$, then $\mathcal{I}_G^{(r)} = {[n] \choose r}$ and ${[n] \choose r}$ is intersecting). One of the difficulties in trying to establish the star property lies in determining a largest star. Our counterexample to Conjecture~\ref{conj} indicates that the problem for trees is more difficult than is hoped.


\section{The counterexample} \label{Results}

Let $x_0 = 0$, $x_1 = 1$ and $x_2 = 2$. For any $k \in \mathbb{N}$, let $y_i = 2+i$ for each $i \in [2k]$, let $z_i = 2k + 2 + i$ for each $i \in [2k]$, and let $T_k$ be the graph whose vertex set is 
\[\{x_0, x_1, x_2\} \cup \{y_i \colon i \in [2k]\} \cup \{z_i \colon i \in [2k]\}\]
and whose edge set is 
\[\{x_0x_1, x_0x_2\} \cup \{x_1y_i \colon i \in [k]\} \cup \{x_2y_{i} \colon i \in [k+1,2k]\} \cup \{y_iz_i \colon i \in [2k]\}.\]
%
%

\begin{theorem} \label{treestars} Let $k \in \mathbb{N}$. \\
(a) The graph $T_k$ is a tree, and the leaves of $T_k$ are $z_1, \dots, z_{2k}$.\\
(b) The largest integer $s$ such that $\mathcal{I}_{T_k}^{(s)}(x_0) \neq \emptyset$ is $2k+1$. \\
(c) If $k \geq 3$, then $|\mathcal{I}_{T_k}^{(r)}(z)| < |\mathcal{I}_{T_k}^{(r)}(x_0)|$ for any leaf $z$ of $T_k$ and any $r \in [5,2k+1]$.
\end{theorem}
\textbf{Proof.} (a) is straightforward. 

Let $G = T_k$. Let $Y = \{y_i \colon i \in [2k]\}$ and $Z = \{z_i \colon i \in [2k]\}$.

We have $\{x_0\} \cup Z \in \mathcal{I}_G^{(2k+1)}(x_0)$. Suppose that $S$ is a set in $\mathcal{I}_G^{(s)}(x_0)$. Then $S \backslash \{x_0\} \in {Y \cup Z \choose s-1}$ and $|(S \backslash \{x_0\}) \cap \{y_i, z_i\}| \leq 1$ for each $i \in [2k]$. Thus $s-1 \leq 2k$, and hence $s \leq 2k+1$. Hence (b).

Suppose $k \geq 3$ and $r \in [5,2k+1]$. Let $\mathcal{J} = \mathcal{I}_G^{(r)}$. Let $\mathcal{E} = \{I \in \mathcal{J} \colon x_0, z_1 \in I\}$. 
Let 
\begin{align} \mathcal{A}_1 &= \{I \in \mathcal{J}(x_0) \colon y_1 \in I\}, \nonumber \\
\mathcal{A}_2 &= \{I \in \mathcal{J}(x_0) \colon y_1, z_1 \notin I\}, \nonumber \\ 
\mathcal{B}_1 &= \{I \in \mathcal{J}(z_1) \colon x_0 \notin I, x_1 \in I, x_2 \notin I\}, \nonumber \\
\mathcal{B}_2 &= \{I \in \mathcal{J}(z_1) \colon x_0 \notin I, x_1 \notin I, x_2 \in I\}, \nonumber \\ 
\mathcal{B}_3 &= \{I \in \mathcal{J}(z_1) \colon x_0 \notin I, x_1, x_2 \in I\}, \nonumber \\
\mathcal{B}_4 &= \{I \in \mathcal{J}(z_1) \colon x_0, x_1, x_2 \notin I\}. \nonumber 
\end{align} 
We have $\mathcal{J}(x_0) = \mathcal{E} \cup \mathcal{A}_1 \cup \mathcal{A}_2$ and $\mathcal{J}(z_1) = \mathcal{E} \cup \mathcal{B}_1 \cup \mathcal{B}_2 \cup \mathcal{B}_3 \cup \mathcal{B}_4$. 
Since $y_1z_1 \in E(G)$, $\{y_1, z_1\} \nsubseteq I$ for each $I \in \mathcal{J}$. Thus $\mathcal{E}$, $\mathcal{A}_1$ and $\mathcal{A}_2$ are pairwise disjoint, and hence 
\begin{equation} |\mathcal{J}(x_0)| = |\mathcal{E}| + |\mathcal{A}_1| + |\mathcal{A}_2|. \label{eqnx0}
\end{equation} 
Since $\mathcal{E}$, $\mathcal{B}_1$, $\mathcal{B}_2$, $\mathcal{B}_3$ and $\mathcal{B}_4$ are pairwise disjoint, 
\begin{equation} |\mathcal{J}(z_1)| = |\mathcal{E}| + |\mathcal{B}_1| + |\mathcal{B}_2| + |\mathcal{B}_3| + |\mathcal{B}_4|. \label{eqnz1}
\end{equation} 

Let $Y' = Y \backslash \{y_1\}$ and $Z' = Z \backslash \{z_1\}$. Since $x_0x_1, x_0x_2 \in E(G)$, we have $\{x_0, x_1\}$, $\{x_0, x_2\} \nsubseteq I$ for each $I \in \mathcal{J}$. Thus $\{A \backslash \{x_0\} \colon A \in \mathcal{A}_2\} = \mathcal{I}_G^{(r-1)} \cap {Y' \cup Z' \choose r-1} = \{B \backslash \{z_1\} \colon B \in \mathcal{B}_4\}$, and hence 
\begin{equation} |\mathcal{A}_2| = |\mathcal{B}_4|. \label{eqn24} 
\end{equation}
%

Let $Y_1 = \{y_i \colon i \in [2,k]\}$ and $Y_2 = \{y_i \colon i \in [k+1,2k]\}$. Let 
\begin{align} \mathcal{A}_1' &= \{A \backslash \{x_0,y_1\} \colon A \in \mathcal{A}_1\}, \nonumber \\
\mathcal{B}_1' &= \{B \backslash \{z_1, x_1\} \colon B \in \mathcal{B}_1\}, \nonumber \\
\mathcal{B}_2' &= \{B \backslash \{z_1, x_2\} \colon B \in \mathcal{B}_2\}, \nonumber \\
\mathcal{B}_3' &= \{B \backslash \{z_1, x_1, x_2\} \colon B \in \mathcal{B}_3\}. \nonumber 
\end{align} 
We have $\mathcal{A}_1' = \mathcal{I}_G^{(r-2)} \cap {Y' \cup Z' \choose r-2}$, $\mathcal{B}_1' = \mathcal{I}_G^{(r-2)} \cap {Y_2 \cup Z' \choose r-2}$, $\mathcal{B}_2' = \mathcal{I}_G^{(r-2)} \cap {Y_1 \cup Z' \choose r-2}$ and $\mathcal{B}_3' = {Z' \choose r-3}$. Let $\mathcal{C} = \{I \in \mathcal{I}_G^{(r-2)} \cap {Y' \cup Z' \choose r-2} \colon I \cap Y_1 \neq \emptyset \neq I \cap Y_2\}$. Thus $\mathcal{A}_1' = \mathcal{B}_1' \cup \mathcal{B}_2' \cup \mathcal{C}$. We have $(\mathcal{B}_1' \cup \mathcal{B}_2') \cap \mathcal{C} = \emptyset$ and $\mathcal{B}_1' \cap \mathcal{B}_2' = {Z' \choose r-2}$. Thus 
\[|\mathcal{A}_1'| = |\mathcal{B}_1' \cup \mathcal{B}_2'| + |\mathcal{C}| = |\mathcal{B}_1'| + |\mathcal{B}_2'| - |\mathcal{B}_1' \cap \mathcal{B}_2'| + |\mathcal{C}| = |\mathcal{B}_1'| + |\mathcal{B}_2'| - {2k-1 \choose r-2} + |\mathcal{C}|. \]
Let $a = |\mathcal{A}_1'| - (|\mathcal{B}_1'| + |\mathcal{B}_2'| + |\mathcal{B}_3'|)$. Then 
\begin{align} a &= |\mathcal{C}| - {2k-1 \choose r-2} - {2k-1 \choose r-3} \nonumber \\
&= \sum_{i=1}^{k-1} \sum_{j=1}^k {|Y_1| \choose i}{|Y_2| \choose j} {|Z'|-i-j \choose r-2-i-j} - {2k-1 \choose r-2} - {2k-1 \choose r-3} \nonumber \\
&= \sum_{i=1}^{k-1} \sum_{j=1}^k {k-1 \choose i}{k \choose j} {2k-1-i-j \choose r-2-i-j} - {2k-1 \choose r-2} - {2k-1 \choose r-3}. \label{eqna} 
\end{align}

We show that $a > 0$. If $r = 2k+1$, then 
\[a = \sum_{i=1}^{k-1} {k-1 \choose i} \sum_{j=1}^k {k \choose j} - 2k = (2^{k-1}-1)(2^k-1) - 2k > 0.\] 
Suppose $r \leq 2k$. We have  
\begin{align} a &\geq {k-1 \choose 1}{k \choose 1} {2k-3 \choose r-4} - {2k-1 \choose r-2} - {2k-1 \choose r-3} \nonumber \\
&= {2k-3 \choose r-4}\left( k(k-1) - \frac{(2k-1)(2k-2)}{(r-2)(r-3)} - \frac{(2k-1)(2k-2)}{(r-3)(2k+2-r)} \right). \nonumber
\end{align}
If $r \geq 6$, then 
\begin{align} a &\geq {2k-3 \choose r-4}\left( k(k-1) - \frac{(2k-1)(2k-2)}{(4)(3)} - \frac{(2k-1)(2k-2)}{(3)(2)} \right) \nonumber \\
&= {2k-3 \choose r-4}\left( k(k-1) - \frac{(2k-1)(k-1)}{2}\right) > 0. \nonumber
\end{align}
If $r=5$ and $k \geq 5$, then 
\begin{equation} a \geq {2k-3 \choose r-4}\left( k(k-1) - \frac{(2k-1)(2k-2)}{(3)(2)} - \frac{(2k-1)(2k-2)}{(2)(7)} \right) > 0. \nonumber
\end{equation}
If $r=5$ and $3 \leq k \leq 4$, then $a>0$ is easily obtained from (\ref{eqna}).
%
%



Since $a > 0$, $|\mathcal{A}_1'| > |\mathcal{B}_1'| + |\mathcal{B}_2'| + |\mathcal{B}_3'|$. Now $|\mathcal{A}_1'| = |\mathcal{A}_1|$, $|\mathcal{B}_1'| = |\mathcal{B}_1|$, $|\mathcal{B}_2'| = |\mathcal{B}_2|$ and $|\mathcal{B}_3'| = |\mathcal{B}_3|$. Thus $|\mathcal{A}_1| > |\mathcal{B}_1| + |\mathcal{B}_2| + |\mathcal{B}_3|$. By (\ref{eqnx0}), (\ref{eqnz1}) and (\ref{eqn24}), it follows that $|\mathcal{J}(x_0)| > |\mathcal{J}(z_1)|$. Clearly, for each $i \in [2k]$, we have $|\mathcal{J}(z_i)| = |\mathcal{J}(z_1)|$, and hence $|\mathcal{J}(z_i)| < |\mathcal{J}(x_0)|$. By (a), (c) follows.~\hfill{$\Box$} \\


If $I$ is a maximal independent set of $T_k$, then $|I \cap \{x_0,x_1,x_2\}| \geq 1$ and $|I \cap \{y_i,z_i\}| = 1$ for each $i \in [2k]$. Thus $\mu(T_k) = 2k+1$. Therefore, if $5 \leq r \leq (2k+1)/2$, 
then the condition $\mu(T_k) \geq 2r$ of the HT Conjecture is satisfied, but, by Theorem~\ref{treestars}, no leaf of $T_k$ yields a star of $\mathcal{I}_{T_k}^{(r)}$ of maximum size.

\end{document}